\documentclass[11pt,reqno]{amsart}
\usepackage{amsthm,amssymb,amsfonts}
\usepackage{mathrsfs}
\usepackage[nosumlimits]{mathtools}
\usepackage[margin=1in]{geometry}
\usepackage{enumerate,tikz,xcolor,microtype,hyperref,needspace}
\hypersetup{pdftitle={Single-graph rational exponents in mixed multicolour Ramsey theory},
 pdfauthor={Yuang Yang},
 pdfsubject={Symmetric random algebraic graphs and rational Ramsey exponents},
 colorlinks,
 linkcolor={red!50!black},
 citecolor={blue!50!black},
 urlcolor={blue!80!black}}
\theoremstyle{plain}
\newtheorem{theorem}{Theorem}[section]
\newtheorem{proposition}[theorem]{Proposition}
\newtheorem{lemma}[theorem]{Lemma}
\newtheorem{corollary}[theorem]{Corollary}
\theoremstyle{definition}
\newtheorem{question}[theorem]{Question}
\DeclareMathOperator{\ex}{ex}
\DeclareMathOperator{\tr}{tr}
\DeclareMathOperator{\tdeg}{trdeg}
\newcommand{\Nind}{N_{\mathrm{ind}}}

\begin{document}
\title[Single-graph rational Ramsey exponents]{Single-graph rational exponents in mixed multicolour Ramsey theory}
\author{Yuang Yang}
\address{School of Mathematical Sciences, University of Science and Technology of China, Hefei 230026, China}
\email{Yangya2023@mail.ustc.edu.cn}
\date{}
\subjclass[2020]{05C55, 05C35, 05D10, 05C80}
\keywords{mixed Ramsey numbers, random algebraic graphs, rooted graphs, rational exponents}
\begin{abstract}
    For every rational number $\gamma>1$, we construct a fixed bipartite graph $H_\gamma$ such that $$r_k(H_\gamma;K_N)=\Theta_{\gamma,k}\bigl((N/\log N)^\gamma\bigr) \qquad\text{for every fixed }k\ge3.$$
    Here $r_k(H;K_N)=r(H,\ldots,H,K_N)$, with $k$ copies of $H$.
    The lower bound uses a symmetric random algebraic graph which avoids a fixed rooted power and has sufficiently strong local density for the Alon--R\"odl method.
    The main algebraic step establishes a nonzero obstruction within the full symmetric coefficient space; the upper bound uses the rooted graphs from the proof of the rational exponents conjecture.
    We also obtain induced even-cycle asymptotics in the same construction.
\end{abstract}
\maketitle

\section{Introduction}

For a graph $H$, let $r_k(H;K_N)=r(H,\ldots,H,K_N)$, with $k$ copies of $H$.
Thus $r_k(H;K_N)$ is the least $n$ such that every $(k+1)$-colouring of $K_n$ contains $H$ in one of the first $k$ colours or $K_N$ in the last colour.
When $H$ is bipartite, the asymptotic behaviour of these numbers can be more accessible for several forbidden colours than for one.
Alon and R\"odl \cite[Theorems~3.4, 3.7--3.9]{alon2005sharp} proved, for every fixed $k\ge3$, estimates of the form \begin{equation}\label{eq:intro-target} r_k(H;K_N)=\Theta\bigl((N/\log N)^\gamma\bigr) \end{equation} for $(H,\gamma)=(C_4,2),(C_6,3/2),(C_{10},5/4)$, and for $(H,\gamma)=(K_{s,t},s)$ when $t\ge s\ge2$ and $t\ge(s-1)!+1$.
Their method uses random relabellings of an $H$-free graph with few large independent sets.
Further results for complete bipartite graphs were obtained by Lenz and Mubayi \cite{lenz2014multicolor}.

These examples suggest asking which exponents can occur in \eqref{eq:intro-target} for a single fixed graph.
We prove that every rational number greater than one occurs.

\begin{theorem}\label{cor:intro-rational}
    For every rational $\gamma>1$, there is a finite bipartite graph $H_\gamma$ such that $$r_k(H_\gamma;K_N)=\Theta_{\gamma,k}\bigl((N/\log N)^\gamma\bigr) \qquad\text{for every fixed }k\ge3.$$
    The graph $H_\gamma$ is independent of $k$, and the constant and threshold in the lower bound can be chosen uniformly over $k\ge3$.
\end{theorem}

Extremal estimates give Ramsey upper bounds and identify candidate forbidden graphs.
By Alon and R\"odl's upper-bound principle \cite[Lemma~3.3]{alon2005sharp}, $$\ex(n,H)=O(n^{2-1/\gamma}) \quad\Longrightarrow\quad r_k(H;K_N)=O_{H,k,\gamma}\bigl((N/\log N)^\gamma\bigr).$$
The proof of the rational exponents conjecture by Conlon, Janzer, Jiang, Kim, Lee, Liu, Longbrake and Yepremyan \cite{conlon2026rational} supplies such extremal bounds through powers of balanced rooted bipartite graphs.
For integers $b>a\ge1$, their construction gives a rooted graph $(F,R)$ of density $b/a$ such that, for every fixed $t\ge1$, $\ex(n,F^{(t)})=O(n^{2-a/b})$.
A rooted power identifies the corresponding roots of $t$ copies and makes no other vertex identifications.
We record the precise input in Section~2.

The extremal bound does not by itself provide the corresponding Ramsey lower bound.
A dense $F^{(t)}$-free graph may have many large independent sets.
Indeed, a bipartite construction has an independent vertex class of linear size.
For the Alon--R\"odl argument we use a construction on a single vertex set, with no prescribed bipartition, in which forbidden-subgraph avoidance is accompanied by local density on all sufficiently large vertex sets.

We use the random algebraic method introduced by Bukh \cite{bukh2015random}, together with the rooted balance framework of Bukh and Conlon \cite{bukh2015rational}.
Their construction gives extremal lower bounds for finite forbidden families, and hence for individual rooted powers.
Symmetric random algebraic constructions with spectral and local-density estimates were developed by Chen, Liu, Yang and Ye \cite{chen2025random}; in particular, they used these estimates to obtain mixed Ramsey bounds for complete bipartite graphs.
Related uses of pseudorandomness and algebraic structure in Ramsey theory appear in \cite{he2020multicolor,sudakov2022ramsey}.

For any given $(a,b)$-balanced rooted graph, the additional task here is to exclude a fixed rooted power in the full symmetric coefficient space while retaining the spectral estimate.
We allow the fixed polynomial degree to depend on the rooted graph and the required trace moment order, so that both interpolation requirements can be met.
An obstruction obtained for unrestricted polynomials in an ordered bipartite model may vanish identically on the symmetric subspace.
Our incidence argument adapts the coefficient-space dimension method of Bukh \cite{bukh2024extremal}; rooted-fibre arguments also occur in \cite{bukh2015rational,conlon2019fewpaths}.
The argument establishes the existence of a fixed power, but gives no explicit bound on $t$.

Choose $a$ independent symmetric forms of sufficiently large fixed bidegree on $\mathbb{P}^b\times\mathbb{P}^b$ over $\mathbb{F}_q$, and join distinct points where all forms vanish.
The affine restriction has $q^b$ vertices and edge probability $q^{-a}$.
Symmetric interpolation gives independence on bounded edge sets, even when their edges share endpoints.
The standard trace method then gives spectral mixing with probability bounded away from zero.

The same interpolation property also controls rooted embeddings.
For a union of copies with common roots, interpolation makes the distinct edge equations independent linear conditions on each coefficient block over the field generated by the vertex coordinates.
Rooted balance then cancels the internal-vertex contribution, including that from additional overlaps, in the transcendence-degree count.
If a fibre were positive-dimensional, sufficiently many independent generic points in it would supply more transcendence than the root coordinates allow.
Thus, in the coefficient-generic graph, the fibre over every geometric root placement is finite.
A uniform bound on these fibres excludes one fixed power over the generic coefficient field.
By Chevalley's theorem, the bad coefficient image is constructible; since it misses the generic point, its Zariski closure is proper.
A nonzero polynomial vanishing on that closure gives a failure probability $O(q^{-1})$, so its nonvanishing event intersects the spectral event.
This is Theorem~\ref{thm:clean-power-spectral-graph}.

The independent-set counting lemma of Kohayakawa, Lee, R\"odl and Samotij \cite{Kohayakawa2015} converts local density into a bound on the number $I_N(G)$ of independent $N$-sets when $N$ is a sufficiently large constant multiple of $q^a\log q$.
The trace moment order can be chosen so that $$I_N(G)^3<\binom{q^b}{N}^{2}.$$
Three random relabellings then yield the lower bound by the Alon--R\"odl method, and the extremal input supplies the upper bound for the same forbidden graph.
The argument for a general balanced rooted graph is stated in Theorem~\ref{H_r}.

For two forbidden colours, the relabelling criterion requires the corresponding logarithmic ratio to be below $1/2$.
Our counting bound approaches $1/2$ from above as the trace moment order increases, whereas three colours require only a ratio below $2/3$.
This is a limitation of the present estimate, and leaves the following question.

\begin{question}\label{q:two-colours}
    Given a rooted graph $(F,R)$ satisfying the hypotheses of Theorem~\ref{H_r}, can one choose a fixed rooted power $H=F^{(t)}$ for which $r_2(H;K_N)=\Theta((N/\log N)^{b/a})$?
\end{question}

After the Ramsey proof, Section~\ref{sec:induced-subgraph-counts} records a further use of finite independence: induced-copy moments, and induced even-cycle asymptotics, agree with those of the independent-edge model in the stated range.
These counting properties can be realized together with forbidden-power avoidance, degree concentration and spectral mixing.

\section{Rooted graphs and Ramsey tools}

All graphs are simple.
An embedding is an injective graph homomorphism; it need not be induced.
We write $e(U)$ for the number of edges induced by $U$, and $e(X,Y)$ for the number of ordered adjacent pairs in $X\times Y$.
Thus edges within $X\cap Y$ are counted twice.
All logarithms are natural, and asymptotic notation refers to the displayed growing parameter with the other parameters fixed.

A rooted graph $(F,R)$ has internal vertex set $A=V(F)\setminus R$, which we require to be nonempty.
For $S\subseteq A$, let $e_F(S)$ count the edges incident with at least one vertex of $S$.
For integers $0<a\le b$, we say that $(F,R)$ is \emph{$(a,b)$-balanced} if \begin{equation}\label{eq:intro-balance} b|S|\le a e_F(S)\qquad(S\subseteq A). \end{equation}
Edges between roots are allowed and do not enter this condition.
The rooted power $F^{(t)}$ consists of $t$ labelled copies of $F$ with corresponding roots identified and internal vertex sets pairwise disjoint.
This is sometimes called the clean rooted power.
We use \emph{rooted power} in this sense throughout, reserving \emph{common-root union} for copies with additional internal overlaps.

\begin{figure}[ht]
    \centering
    \begin{tikzpicture}[x=.8cm,y=.7cm,
        root/.style={circle,fill=black,inner sep=1.8pt},
        internal/.style={circle,draw,fill=white,inner sep=1.8pt},
        every label/.style={font=\small}]
        \begin{scope}
            \node[root,label=left:$r_1$] (a) at (0,0) {};
            \node[internal,label=above:$u$] (b) at (1,0) {};
            \node[internal,label=above:$v$] (c) at (2,0) {};
            \node[root,label=right:$r_2$] (d) at (3,0) {};
            \draw (a)--(b)--(c)--(d);
            \node at (1.5,-1.4) {\small $F$};
        \end{scope}
        \begin{scope}[xshift=4.3cm]
            \node[root] (a) at (0,0) {};
            \node[root] (d) at (3,0) {};
            \node[internal] (b) at (1,.7) {};
            \node[internal] (c) at (2,.7) {};
            \node[internal] (e) at (1,-.7) {};
            \node[internal] (f) at (2,-.7) {};
            \draw (a)--(b)--(c)--(d)--(f)--(e)--(a);
            \node at (1.5,-1.4) {\small $F^{(2)}\cong C_6$};
        \end{scope}
        \begin{scope}[xshift=8.5cm]
            \node[root] (a) at (0,0) {};
            \node[internal] (b) at (1,0) {};
            \node[internal] (c) at (2,.7) {};
            \node[internal] (e) at (2,-.7) {};
            \node[root] (d) at (3,0) {};
            \draw (a)--(b)--(c)--(d)--(e)--(b);
            \node at (1.5,-1.4) {\small Additional overlap};
        \end{scope}
    \end{tikzpicture}
    \caption{Filled vertices are roots. A rooted power identifies only corresponding roots; a common-root union may have further overlaps.}
    \label{fig:clean-rooted-path}
\end{figure}
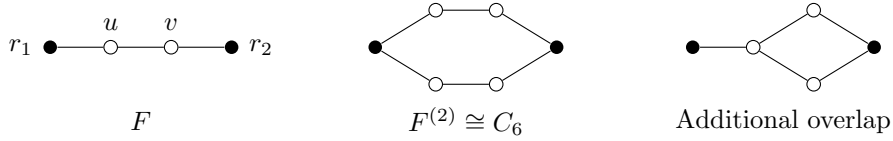

For example, a path $r_1uvr_2$ rooted at its endpoints has $F^{(2)}\cong C_6$.
If the two copies also identify $u_1=u_2$, their union is a four-cycle with a pendant edge, as shown in Figure~\ref{fig:clean-rooted-path}.
Such unions belong to the power family considered by Bukh and Conlon \cite{bukh2015rational}, but need not be rooted powers in our notation.

We use the following consequence of the rooted-graph construction in \cite[Section~2]{conlon2026rational}.

\begin{theorem}[Conlon et al.]\label{thm:clean-power-input}
    For every pair of integers $b>a\ge1$, there is a finite rooted bipartite graph $(F,R)$, with $V(F)\setminus R\ne\varnothing$, satisfying \eqref{eq:intro-balance} and, for every fixed $t\ge1$, \begin{equation}\label{eq:intro-clean-upper} \ex(n,F^{(t)})=O_{F,t}(n^{2-a/b}). \end{equation}
\end{theorem}

In the cited paper, Proposition~2.3 supplies a balanced graph of density $b/a$, and equation~(4) in Section~2.4 gives \eqref{eq:intro-clean-upper}.
Its definition of a rooted power agrees with ours.

We also record the three Ramsey tools used in Section~5.
Write $I_s(G)$ for the number of independent $s$-sets in $G$.

\begin{lemma}[Alon--R\"odl {\cite[Lemma~3.1]{alon2005sharp}}]\label{indep-set}
    Let $G$ be an $n$-vertex graph and $k\ge2$.
    If $$I_s(G)^k<\binom ns^{k-1},$$ then $k$ relabellings of $G$ on the same vertex set have no common independent $s$-set.
    If $G$ is $H$-free, then $r_k(H;K_s)>n$.
\end{lemma}

\begin{proof}
    For independent uniform relabellings, the expected number of common independent $s$-sets is $I_s(G)^k/\binom ns^{k-1}<1$.
    Assign each edge present in their union to one relabelling containing it, and give the remaining edges the last colour.
\end{proof}

\Needspace{10\baselineskip}
\begin{lemma}[Kohayakawa--Lee--R\"odl--Samotij {\cite[Lemma~3.1]{Kohayakawa2015}}]\label{container}
    Let $G$ have $n$ vertices, let $f,M$ be positive integers, and let $0<\beta<1$.
    Suppose that $ne^{-\beta f}<M$ and $2e(U)\ge\beta|U|^2$ whenever $|U|\ge M$.
    Then, for every integer $s\ge f$, $$I_s(G)\le\binom nf\binom M{s-f}.$$
    Here a binomial coefficient is zero if its lower argument exceeds its upper argument.
\end{lemma}

This is the cited lemma with $q=f$ and $r=s-f$; our density hypothesis implies its hypothesis $e(U)\ge\beta\binom{|U|}{2}$.

\begin{lemma}[Alon--R\"odl {\cite[Lemma~3.3]{alon2005sharp}}]\label{upper-bound}
    Let $H$ be a fixed bipartite graph and $\gamma>1$.
    If $\ex(n,H)=O(n^{2-1/\gamma})$, then, for every fixed $k\ge1$, $$r_k(H;K_N)=O_{H,k,\gamma}\bigl((N/\log N)^\gamma\bigr).$$
\end{lemma}

\section{Finite edge independence and spectral mixing}

For a field $\kappa$, put $$R_m(\kappa)=\kappa[Z_0,\ldots,Z_b]_m,\qquad \mathscr{S}_m(\kappa)=(R_m(\kappa)\otimes_\kappa R_m(\kappa))^\tau, \qquad \tau(P)(X,Y)=P(Y,X).$$
We use the full swap-invariant space, also in characteristic two.
If $M_1,\ldots,M_d$ are the degree-$m$ monomials, it has basis $$M_i(X)M_i(Y)\quad(1\le i\le d),\qquad M_i(X)M_j(Y)+M_j(X)M_i(Y)\quad(i<j).$$
In particular, this space commutes with field extension.
No division by two is involved.
Evaluation at a projective point uses a fixed nonzero representative; vanishing is independent of the choice.

Choose $a$ independent uniform members $P_1,\ldots,P_a$ of $\mathscr{S}_m(\mathbb{F}_q)$.
Let $\widetilde G_q$ be the graph on $\mathbb{P}^b(\mathbb{F}_q)$ joining distinct points $x,y$ when $P_j(x,y)=0$ for every $j$, and let $G_q$ be its restriction to the affine chart $x\mapsto[1:x]$.
The edge probability is $p=q^{-a}$.
We first establish the interpolation properties used in both the trace estimate and the algebraic argument.

\begin{lemma}[Symmetric interpolation]\label{edge-pr-sym}\label{lem:geometric-symmetric-interpolation}
    Let $\kappa$ be any field, and let $e_1,\ldots,e_s$ be distinct unordered pairs of distinct points of $\mathbb{P}^b(\kappa)$.
    If $m\ge s-1$, the evaluation map $$\mathscr{S}_m(\kappa)\longrightarrow\kappa^s,\qquad P\longmapsto(P(e_1),\ldots,P(e_s))$$ is surjective.
    Consequently, in the finite-field model any $s$ distinct edge indicators are independent Bernoulli variables of mean $q^{-a}$, whether or not their endpoints overlap.
\end{lemma}

\begin{proof}
    Fix $e_i=\{u_i,v_i\}$.
    For each $j\ne i$, choose $z_{ij}\in e_j\setminus\{u_i,v_i\}$.
    In the space of linear forms vanishing at $z_{ij}$, the conditions of vanishing at $u_i$ and at $v_i$ define two proper subspaces.
    Since two proper subspaces cannot cover a vector space, there is a form $\ell_{ij}$ vanishing at $z_{ij}$ and nonzero at both endpoints of $e_i$.
    Likewise choose $M_i$ nonzero at both endpoints, and set $$h_i=\left(\prod_{j\ne i}\ell_{ij}\right)M_i^{m-s+1}, \qquad H_i(X,Y)=h_i(X)h_i(Y).$$
    Then $H_i$ is symmetric, is nonzero on $e_i$, and vanishes on every other $e_j$.
    Rescaling gives the coordinate vectors in the image.
    This construction works in characteristic two because $h_i\otimes h_i$ belongs to the full invariant space.
    A uniform vector has uniform image under a surjective linear map; apply this separately to the $a$ coefficient blocks.
\end{proof}

For induced counts we need all pairs on a bounded vertex set.
The next lemma achieves this with a degree bound in terms of the number of vertices, rather than the number of pairs.

\begin{lemma}[Vertex-supported independence]\label{vertex-supported-independence}
    Let $x_1,\ldots,x_r$ be distinct points of $\mathbb{P}^b(\mathbb{F}_q)$.
    If $m\ge r-1$, evaluation on their $\binom r2$ unordered pairs maps $\mathscr{S}_m(\mathbb{F}_q)$ surjectively onto $\mathbb{F}_q^{\binom r2}$.
    All edge indicators on these vertices are therefore mutually independent.
\end{lemma}

\begin{proof}
    Products of separating linear forms, padded to degree $m$ and rescaled, give $h_i\in R_m(\mathbb{F}_q)$ with $h_i(x_j)=\delta_{ij}$.
    For $i<j$, $h_i(X)h_j(Y)+h_j(X)h_i(Y)$ evaluates to the coordinate vector for $\{x_i,x_j\}$, also in characteristic two.
\end{proof}

The following trace argument is standard for graphs with finite edge independence; see \cite{alon2008kwise} and, in the symmetric algebraic setting, \cite[Appendix~A]{chen2025random}.

\begin{theorem}[Spectral mixing]\label{thm:spectral-mixing-rag}
    Fix integers $b>a\ge1$, $h\ge1$, and $m\ge2h-1$.
    For $\widetilde G_q$, put $n=|\mathbb{P}^b(\mathbb{F}_q)|$, $p=q^{-a}$, and $B=A(\widetilde G_q)-p(J-I)$.
    There is $K=K(a,b,h)$ such that, for every sufficiently large $q$, with probability at least $1/2$, \begin{equation}\label{eq:projective-spectral} \|B\|\le K n^{1/(2h)}\sqrt{np}, \qquad \|B\|\le Kq^{(b-a)/2+b/(2h)}. \end{equation}
    On this event, simultaneously for all $X,Y$, \begin{equation}\label{eq:mixing} \left|e(X,Y)-p\bigl(|X||Y|-|X\cap Y|\bigr)\right| \le\|B\|\sqrt{|X||Y|}. \end{equation}
    In particular, for disjoint $X,Y$ and $\varepsilon>0$, $e(X,Y)=(1\pm\varepsilon)p|X||Y|$ whenever $\sqrt{|X||Y|}\ge(K/\varepsilon)q^{(b+a)/2+b/(2h)}$.
\end{theorem}

\begin{proof}
    Write $\xi_e$ for the edge indicators.
    Since $B$ is symmetric with zero diagonal, $$\|B\|^{2h}\le\tr(B^{2h}) =\sum_{v_0,\ldots,v_{2h-1}} \prod_{i=0}^{2h-1}B_{v_i v_{i+1}},\qquad v_{2h}=v_0.$$
    Terms with consecutive equal vertices vanish.
    Each remaining term is a closed walk of length $2h$, involving at most $2h$ distinct unordered pairs.
    Lemma~\ref{edge-pr-sym} makes their indicators independent.
    A pair traversed exactly once contributes the factor $\mathbb{E}(\xi_e-p)=0$, so every edge of a contributing walk is traversed at least twice.

    If such a walk has $s$ distinct edges and $r$ vertices, its support is connected and $r\le s+1\le h+1$.
    For a fixed $s$, there are at most $$\sum_{r=1}^{s+1}\binom nr r^{2h}\le C_h n^{s+1}$$ walks.
    Moreover, $\mathbb{E}|\xi_e-p|^j\le2p$ for every $j\ge2$, so a walk with $s$ distinct edges contributes at most $(2p)^s$ in absolute value.
    Consequently, since $np\to\infty$, \begin{equation}\label{eq:trace-bound} \mathbb{E}\tr(B^{2h}) \le C_h\sum_{s=1}^h n^{s+1}(2p)^s \le C_h' n(np)^h. \end{equation}
    Choosing $K_0^{2h}\ge2C_h'$, Markov's inequality gives $$\Pr\!\left(\|B\|>K_0n^{1/(2h)}\sqrt{np}\right) \le\frac{\mathbb{E}\tr(B^{2h})}{K_0^{2h}n(np)^h}\le\frac12.$$
    Since $n=(1+o(1))q^b$, increasing the constant gives \eqref{eq:projective-spectral}.
    Finally, the left side of \eqref{eq:mixing} is $\left|\mathbf{1}_X^{\mathsf{T}}B\mathbf{1}_Y\right|$.
    The operator norm inequality proves \eqref{eq:mixing} and its stated consequence.
\end{proof}

\begin{corollary}[Local density]\label{cor:local-density-from-mixing}\label{cor:jumbledness-interface}
    On the spectral event, there is $C_0=C_0(a,b,h)$ such that $$|U|\ge C_0q^{(b+a)/2+b/(2h)} \quad\Longrightarrow\quad 2e(U)\ge\frac p4|U|^2 \quad\text{and}\quad e(U)\ge\frac p4\binom{|U|}{2}.$$
    Also, $\left|e(X,Y)-p|X||Y|\right|\le\beta'\sqrt{|X||Y|}$ with $\beta'=\|B\|+p=O(q^{(b-a)/2+b/(2h)})$.
    In particular, $\beta'=o(p^\lambda n)$ for every fixed $\lambda<(b+a)/(2a)-b/(2ah)$.
\end{corollary}

\begin{proof}
    For $u=|U|$, \eqref{eq:mixing} gives $2e(U)\ge p(u^2-u)-\|B\|u$, which is at least $pu^2/4$ when $u\ge2$ and $u\ge4\|B\|/p$.
    Absorbing the diagonal term in \eqref{eq:mixing} gives $\beta'=\|B\|+p$; comparing powers of $q$ proves the last assertion.
\end{proof}

We will also use degree concentration in the simultaneous realization in Section~6.
It is not needed for the Ramsey lower bound.

\begin{proposition}[Uniform degrees]\label{prop:uniform-degrees}
    Fix $b>a\ge1$, an even integer $r>2b/(b-a)$, and $m\ge r-1$.
    Write $n$ for the size of the projective vertex set or its affine chart.
    On either vertex set, with probability $1-o(1)$, $$\deg_G(v)=(1+o(1))pn\qquad\text{uniformly in }v.$$
    If $m\ge2h-1$ as well, this event and the spectral event occur together with probability at least $1/2-o(1)$.
\end{proposition}

\begin{proof}
    Put $\mu=(n-1)p$.
    The $r$-wise independence of the edges at a fixed vertex makes every term with a singleton index vanish in the $r$-th centred degree moment.
    Terms with $j\le r/2$ distinct indices contribute $O_r(\mu^j)$ in total, giving $\mathbb{E}|\deg_G(v)-\mu|^r\le C_r\mu^{r/2}$.
    Set $\xi=r(b-a)/2-b>0$ and $\varepsilon_q=q^{-\xi/(2r)}$.
    A union bound and Markov's inequality give $$\Pr\bigl(\exists v:\,|\deg_G(v)-\mu|>\varepsilon_q\mu\bigr) \le C_r n\varepsilon_q^{-r}\mu^{-r/2} =O(q^{-\xi/2})=o(1).$$
    This argument applies directly to either vertex set.
    Intersecting with an event of probability at least $1/2$ proves the last assertion.
\end{proof}

\section{A symmetric obstruction for rooted powers}\label{sec:symmetric-obstruction}

The combinatorial input is the balance inequality for common-root unions.
The algebraic input is Lemma~\ref{edge-pr-sym}, which makes the distinct image edges impose independent conditions on symmetric coefficients.

\subsection{Balance under overlaps}

\begin{lemma}[Common-root density]\label{lem:common-root-union}
    Let $(F,R)$ be $(a,b)$-balanced.
    Suppose that $\phi_1,\ldots,\phi_s$ are embeddings of $F$ in a simple graph which agree on the labelled roots.
    If $I=\bigcup_i\phi_i(A)$ and $E=\bigcup_i\phi_i(E(F))$, then \begin{equation}\label{eq:common-root-union} b|I|\le a|E|. \end{equation}
    Here $E$ consists of distinct image edges, and arbitrary additional internal overlaps are allowed.
\end{lemma}

\begin{proof}
    Add the copies in order.
    Let $I_{i-1},E_{i-1}$ be the internal vertices and edges already present, and set $S_i=\{x\in A:\phi_i(x)\notin I_{i-1}\}$.
    An internal vertex cannot coincide with a common root image, by injectivity.
    Hence $\phi_i(S_i)$ is exactly the set of new vertices.
    Every edge incident with $S_i$ has a new endpoint and therefore a new image, and distinct edges in one copy have distinct images.
    Thus $$|I_i\setminus I_{i-1}|=|S_i|, \qquad |E_i\setminus E_{i-1}|\ge e_F(S_i).$$
    Apply balance to each $S_i$ and sum.
    Edges between roots are not counted in $e_F(S_i)$, so their repetition causes no difficulty.
\end{proof}

This is the common-root argument of Bukh and Conlon \cite[Lemma~2.2]{bukh2015rational}, stated for the balance condition we use.
Its role is to allow overlap in the embeddings selected below: each additional internal coordinate must still be accounted for by sufficiently many distinct edge equations.

\subsection{Generic rooted fibres}

Fix a field $\kappa$, a basis $B_1,\ldots,B_D$ of $\mathscr{S}_m(\kappa)$, and independent indeterminates $c_{j\nu}$, $1\le j\le a$, $1\le\nu\le D$.
Put $$\mathbb{K}=\kappa(c_{j\nu}),\qquad P_j=\sum_{\nu=1}^D c_{j\nu}B_\nu.$$
For any field containing the coefficients, let $\Gamma_{\mathbf{P}}$ be the graph on its affine $b$-space defined by \begin{equation}\label{eq:symmetric-affine-graph} x\sim y\quad\Longleftrightarrow\quad x\ne y\ \text{ and }\ P_j([1:x],[1:y])=0\quad(1\le j\le a). \end{equation}

\begin{proposition}[Generic fibre finiteness]\label{prop:generic-rooted-fibres}
    Let $(F,R)$ be $(a,b)$-balanced, put $L_0=b|R|+1$, and suppose \begin{equation}\label{eq:generic-degree-budget} m\ge L_0e(F)-1. \end{equation}
    For the coefficient-generic graph $\Gamma_{\mathbf{P}}$, every geometric fibre of the map from labelled embeddings of $F$ to their root placements is finite.
    This holds for root placements over every algebraically closed extension of $\mathbb{K}$.
\end{proposition}

\begin{proof}
    Fix a root placement $\rho:R\to\Omega^b$ over an algebraically closed extension $\Omega/\mathbb{K}$.
    A noninjective placement has empty embedding fibre.
    Suppose an injective placement has an infinite fibre $X$.
    The edge equations and distinctness conditions make $X$ constructible over $\mathbb{K}(\rho)$, with coordinates indexed by the internal vertices.

    Over $\overline{\mathbb{K}(\rho)}$, the closure of $X$ has a positive-dimensional irreducible component on which $X$ contains a nonempty open subset.
    A generic point of the $L_0$-fold product of this component gives $L_0$ independent generic points in that open subset, in a common extension field.
    These give embeddings $\phi_1,\ldots,\phi_{L_0}$ and an internal coordinate $z_i$ of each such that $z_1,\ldots,z_{L_0}$ are algebraically independent over $\mathbb{K}(\rho)$.
    Indeed, successive generic points have positive transcendence degree over the field generated by the preceding ones, so at each step some coordinate is transcendental.
    Together with the generic coefficients, the selected coordinates therefore have transcendence degree $aD+L_0$ over $\kappa$.

    Let $I$ and $E$ be the distinct internal vertices and image edges in these embeddings, and put $U=I\cup\rho(R)$.
    By Lemma~\ref{lem:common-root-union} and \eqref{eq:generic-degree-budget}, \begin{equation}\label{eq:generic-union-balance} b|I|\le a|E|,\qquad |E|\le L_0e(F)\le m+1. \end{equation}
    Let $L$ be the field generated over $\kappa$ by the affine coordinates of $U$.
    For each $j$, the edge equations are $$\sum_{\nu=1}^D B_\nu([1:x],[1:y])c_{j\nu}=0 \qquad(\{x,y\}\in E).$$
    The basis $B_\nu$ remains a basis after extension to $L$.
    Lemma~\ref{edge-pr-sym} says that this system has rank $|E|$.
    Consequently $|E|$ coefficients in each block can be expressed over $L$ in terms of the remaining ones.
    We obtain \begin{equation}\label{eq:generic-trdeg-count} \begin{aligned} aD+L_0 &\le \tdeg_\kappa L(c_{j\nu})\\ &\le b|U|+a(D-|E|)\\ &\le aD+b|R|. \end{aligned} \end{equation} The middle line counts at most $b|U|$ vertex coordinates and $a(D-|E|)$ free coefficients.
    Balance cancels the internal-vertex contribution, leaving only the $b|R|$ root coordinates.
    This contradicts $L_0=b|R|+1$.
\end{proof}

The quantifier over root placements is essential: the roots may depend algebraically or transcendently on the generic coefficients.
The argument above allows this by taking the independent internal coordinates over $\mathbb{K}(\rho)$, rather than only over $\mathbb{K}$.

\begin{corollary}\label{cor:generic-clean-power}
    Under the hypotheses of Proposition~\ref{prop:generic-rooted-fibres}, there is a fixed integer $t\ge1$ such that $\Gamma_{\mathbf{P}}(\overline{\mathbb{K}})$ is $F^{(t)}$-free.
\end{corollary}

\begin{proof}
    The embedding scheme is locally closed in affine space over $\mathbb{K}$.
    Its projection to root placements is of finite type and has finite geometric fibres, so it is quasi-finite.
    Its source is quasi-compact; hence the sizes of all geometric fibres have a common bound $M$ \cite[Tag~03JA]{stacks-finite-fibres}.
    An embedding of $F^{(M+1)}$ would give $M+1$ distinct embeddings of $F$ over the same roots.
\end{proof}

\subsection{From generic coefficients to a polynomial obstruction}

Generic avoidance alone is not a probability estimate over $\mathbb{F}_q$.
We need a fixed polynomial whose nonvanishing certifies avoidance after specialization.
Let $\mathcal{C}_m=\mathbb{A}_\kappa^{aD}$ be the symmetric coefficient space.

\begin{proposition}[Symmetric coefficient obstruction]\label{prop:symmetric-coefficient-obstruction}
    Let $(F,R)$ be $(a,b)$-balanced and $m\ge(b|R|+1)e(F)-1$.
    There exist an integer $t\ge1$ and a nonzero polynomial $Q\in\kappa[\mathcal{C}_m]$ such that, for every field extension $L/\kappa$ and every $c\in\mathcal{C}_m(L)$, \begin{equation}\label{eq:principal-obstruction} Q(c)\ne0\quad\Longrightarrow\quad \Gamma_c(\overline L)\text{ is }F^{(t)}\text{-free}. \end{equation}
    With $\kappa$ and the coefficient basis fixed, $t$ and $Q$ depend only on $F,a,b,m$.
\end{proposition}

\begin{proof}
    Take $t$ from Corollary~\ref{cor:generic-clean-power} and write $J=F^{(t)}$.
    In $\mathcal{C}_m\times(\mathbb{A}^b)^{V(J)}$, impose the edge equations for $J$ and remove all vertex diagonals.
    This defines a locally closed incidence scheme $\mathcal{X}_J$.
    Its projection $\pi:\mathcal{X}_J\to\mathcal{C}_m$ is of finite type between Noetherian schemes, so its image $\mathcal{B}_J$ is constructible by Chevalley's theorem \cite[Tag~054K]{stacks-chevalley}.

    The generic point of $\mathcal{C}_m$ does not belong to $\mathcal{B}_J$: otherwise its finite-type fibre would have a point over $\overline{\mathbb{K}}$, contradicting Corollary~\ref{cor:generic-clean-power}.
    A dense constructible subset of the irreducible space $\mathcal{C}_m$ contains a nonempty open subset and hence its generic point.
    Therefore $Z=\overline{\mathcal{B}_J}$ is proper.
    Choose $0\ne Q\in I(Z)$.

    The inverse image of $\mathcal{C}_m\setminus Z$ is the empty scheme, and remains empty after every base change.
    Since the principal open $D(Q)$ is contained in this complement, the embedding fibre at any $L$-point with $Q(c)\ne0$ is empty even after extension to $\overline L$.
    This proves \eqref{eq:principal-obstruction}.
\end{proof}

\subsection{Simultaneous forbidden-power avoidance and mixing}

We now fix $\kappa=\mathbb{F}_2$ and let $q=2^s$.
The full invariant space used in Section~3 makes the interpolation argument valid in this characteristic.
The fixed ratio between successive field sizes will also be useful in the Ramsey argument.

\begin{theorem}[Simultaneous construction]\label{thm:clean-power-spectral-graph}
    Fix integers $b>a\ge1$ and $h\ge1$, and an $(a,b)$-balanced rooted graph $(F,R)$.
    Choose a fixed degree \begin{equation}\label{eq:combined-degree-choice} m\ge\max\{2h-1,(b|R|+1)e(F)-1\}. \end{equation}
    There are an integer $t\ge1$ and constants $K,C_0>0$ such that, for every sufficiently large $q=2^s$, a graph $G_q$ on $n_q=q^b$ vertices simultaneously satisfies:
    \begin{enumerate}[(i)]
        \item $G_q$ is $F^{(t)}$-free;
        \item for $p=q^{-a}$, \begin{equation}\label{eq:affine-spectral-bound} \|A(G_q)-p(J-I)\|\le Kq^{(b-a)/2+b/(2h)}; \end{equation}
        \item every $U\subseteq V(G_q)$ of size at least $C_0q^\theta$, where $\theta=(b+a)/2+b/(2h)$, satisfies \begin{equation}\label{eq:clean-power-local-density} 2e(G_q[U])\ge\frac{q^{-a}}8|U|^2. \end{equation}
    \end{enumerate}
    The same conclusion holds for $q=\ell^s$ with any fixed prime $\ell$, allowing the fixed data to depend on $\ell$.
\end{theorem}

The choices occur in the order $$(F,a,b,h)\longmapsto m\longmapsto(t,Q,K,C_0) \longmapsto q_0\longmapsto q=2^s\ge q_0.$$
In particular, the forbidden power is fixed before $q$ varies.

\begin{proof}
    Proposition~\ref{prop:symmetric-coefficient-obstruction} over $\mathbb{F}_2$ gives $t$ and $Q\ne0$.
    Sample the $a$ forms uniformly over $\mathbb{F}_q$.
    The spectral event $\mathcal{E}_{\mathrm{sp}}$ for $\widetilde G_q$ has probability at least $1/2$.
    The polynomial $Q$ remains nonzero over $\mathbb{F}_q$, and Schwartz--Zippel gives $$\Pr(Q(c)=0)\le\frac{\deg Q}{q},\qquad \Pr(\mathcal{E}_{\mathrm{sp}}\cap\{Q(c)\ne0\}) \ge\frac12-\frac{\deg Q}{q}>0$$ for large $q$.
    On this intersection the affine graph is $F^{(t)}$-free.
    Its centred adjacency matrix is a principal submatrix of the projective one, so it satisfies \eqref{eq:affine-spectral-bound}.
    The local-density conclusion follows from Corollary~\ref{cor:local-density-from-mixing}.
    For a different fixed prime, apply the same argument over $\mathbb{F}_\ell$.
\end{proof}

\section{Single-graph rational Ramsey exponents}\label{sec:single-graph-ramsey}

We first state the general form of the Ramsey argument.

\begin{theorem}\label{H_r}
    Let $b>a\ge1$ be integers, and let $(F,R)$ be a finite $(a,b)$-balanced rooted bipartite graph with nonempty internal vertex set.
    Suppose that $\ex(n,F^{(u)})=O_{F,u}(n^{2-a/b})$ for every fixed $u\ge1$.
    Then there is $t=t(F,a,b)$ such that, with $H=F^{(t)}$, $$r_k(H;K_N)=\Theta_{F,a,b,k}\bigl((N/\log N)^{b/a}\bigr) \qquad\text{for every fixed }k\ge3.$$
    The lower-bound constant and threshold are uniform over $k\ge3$.
    The lower bound uses only balance; the extremal hypothesis is needed only for the value of $t$ selected by the construction.
\end{theorem}

\begin{proof}
    We aim to construct an $H$-free graph $G_q$ on $q^b$ vertices such that \begin{equation}\label{eq:ramsey-goal} I_{N_q}(G_q)^3<\binom{q^b}{N_q}^{2}, \qquad N_q=\Theta(q^a\log q). \end{equation}
    Lemma~\ref{indep-set} will then give the lower bound using three relabellings.

    Choose an integer $h>3b/(b-a)$, then $m$ as in \eqref{eq:combined-degree-choice}.
    Theorem~\ref{thm:clean-power-spectral-graph} gives a fixed $t$ and $F^{(t)}$-free graphs $G_q$ for all sufficiently large powers of two.
    Put \begin{equation}\label{eq:ramsey-delta} n_q=q^b,\qquad \theta=\frac{b+a}{2}+\frac{b}{2h},\qquad \delta_h=\frac{\theta-a}{b-a} =\frac12+\frac{b}{2h(b-a)}<\frac23. \end{equation}
    Thus $a<\theta<b$.
    The point of the choice of $h$ is the last inequality, which leaves room in the entropy comparison for the independent-set encoding.

    With $C_0$ from the construction theorem, set $$M_q=\lceil C_0q^\theta\rceil,\qquad \beta=\frac1{8q^a},\qquad f_q=\lceil c_0q^a\log q\rceil,$$ where $c_0>8(b-\theta)$.
    The local-density hypothesis of Lemma~\ref{container} holds, and $n_qe^{-\beta f_q}\le q^{b-c_0/8}<M_q$ for large $q$.
    Choose $\eta>0$ with $\delta_h+2\eta<2/3$, then choose $c_1>c_0/\eta$ and put $N_q=\lceil c_1q^a\log q\rceil$.
    For large $q$, $$f_q/N_q\le\eta,\qquad 1\le N_q-f_q<M_q,\qquad N_q<n_q.$$
    Writing $I_q=I_{N_q}(G_q)$, we obtain \begin{equation}\label{eq:single-graph-container} I_q\le\binom{n_q}{f_q}\binom{M_q}{N_q-f_q}. \end{equation}

    If $I_q=0$, \eqref{eq:ramsey-goal} is immediate.
    Otherwise, standard binomial estimates give $$\begin{aligned} \log\binom{n_q}{f_q} &\le f_q(b-a+o(1))\log q,\\ \log\binom{M_q}{N_q-f_q} &\le (N_q-f_q)(\theta-a+o(1))\log q,\\ \log\binom{n_q}{N_q} &\ge N_q(b-a+o(1))\log q. \end{aligned}$$ Here $f_q,N_q,N_q-f_q=q^{a+o(1)}$ and $M_q=q^{\theta+o(1)}$.
    Dividing the logarithm of \eqref{eq:single-graph-container} by the last quantity gives \begin{equation}\label{eq:single-graph-entropy-gap} \frac{\log I_q}{\log\binom{n_q}{N_q}} \le \frac{f_q}{N_q} +\left(1-\frac{f_q}{N_q}\right)\delta_h+o(1) \le\delta_h+\eta+o(1)<\frac23. \end{equation}
    This proves \eqref{eq:ramsey-goal}, and hence $r_3(F^{(t)};K_{N_q})>q^b$.

    For any sufficiently large $N$, choose the largest power of two $q$ with $N_q\le N$.
    Since $N_q\le N<N_{2q}$, $$N=\Theta(q^a\log q),\qquad q^b=\Theta\bigl((N/\log N)^{b/a}\bigr).$$
    Monotonicity in $N$, and then in $k$, proves the required lower bound uniformly over $k\ge3$.

    Finally, $H=F^{(t)}$ is bipartite: the bipartition of each copy of $F$ is consistent at corresponding roots.
    The assumed extremal estimate for this $t$, together with Lemma~\ref{upper-bound}, gives the matching upper bound for every fixed $k$.
\end{proof}

\begin{proof}[Proof of Theorem~\ref{cor:intro-rational}]
    Write $\gamma=b/a>1$ in lowest terms.
    Theorem~\ref{thm:clean-power-input} supplies a rooted graph satisfying the hypotheses of Theorem~\ref{H_r}.
    Its power selected by that theorem is the required $H_\gamma$.
\end{proof}

\section{Further applications: induced subgraph counts}\label{sec:induced-subgraph-counts}

Finite independence also determines induced-copy moments on bounded vertex sets.
For each fixed even cycle, a second-moment argument shows that its count in our model is $(1+o(1))$ times its expectation with high probability.
These counting properties can be imposed together with rooted-power avoidance, degree concentration and spectral mixing on the same affine graph, by selecting one coefficient tuple and making no further vertex or edge deletions.
We also record a deterministic consequence of spectral mixing via an induced-cycle supersaturation theorem.

Let $\Gamma_p$ be the independent $p$-random spanning subgraph of a fixed $n$-vertex host $\Gamma$.
For a fixed graph $J$, let $X_J(G)$ count injections $\phi:V(J)\to V(\Gamma)$ such that, for all distinct $u,v\in V(J)$, $$uv\in E(J)\quad\Longleftrightarrow\quad \phi(u)\phi(v)\in E(G).$$
Thus the nonedge requirements concern $G$, not the host $\Gamma$.

\begin{proposition}[Induced-copy moments]\label{prop:induced-moment-transfer}
    Let $G$ be a random spanning subgraph of the $n$-vertex host $\Gamma$.
    Suppose that on every set of at most $L$ vertices the eligible edge indicators are independent Bernoulli variables of mean $p$.
    Put $v=|V(J)|\ge1$.
    For every positive integer $r$ with $rv\le L$, $$\mathbb{E} X_J(G)^r=\mathbb{E} X_J(\Gamma_p)^r.$$
    For $\Gamma=K_n$ and $v\le L$, the common mean is $$\mu_J=(n)_v p^{e(J)}(1-p)^{\binom v2-e(J)}.$$
    If $2v\le L$, the means and variances agree.
    In particular, writing $\mu_J=\mathbb{E} X_J(\Gamma_p)>0$ on any host, $\operatorname{Var}(X_J(\Gamma_p))=o(\mu_J^2)$ implies $X_J(G)=(1+o(1))\mu_J$ with probability $1-o(1)$.
\end{proposition}

\begin{proof}
    Expand $X_J^r$ into products of $r$ induced-copy indicators.
    Each product depends on pairs within at most $rv$ vertices, where the joint law agrees with that of $\Gamma_p$; pairs outside the host are absent in both models.
    Sum the equal expectations.
    The last assertion follows from the cases $r=1,2$ and Chebyshev's inequality.
\end{proof}

More generally, a polynomial model on any host satisfies this hypothesis whenever evaluation from its coefficient space onto all eligible pairs in each such vertex set is surjective.
Choosing $a$ independent uniform polynomials and requiring simultaneous vanishing then gives $p=q^{-a}$.
In our symmetric model, Lemma~\ref{vertex-supported-independence} permits $L=m+1$; in particular, the first two moments for a $v$-vertex graph agree with those of $G(n,p)$ if $m\ge2v-1$.
These identities concern the unconditioned model.
A counting event of probability $1-o(1)$ remains such after conditioning on an event of probability bounded away from zero, although exact moment identities need not survive that conditioning.

Write $\Nind(C_{2\ell},G)$ for the number of unlabelled induced copies of $C_{2\ell}$.

\begin{theorem}[Induced even cycles]\label{thm:sharp-induced-cycles}
    Fix an integer $\ell\ge2$, and let $G=G_n$ be a random graph on $n$ vertices, with $p=p(n)$.
    Suppose that on every set of at most $4\ell$ vertices all edge indicators are independent Bernoulli variables of mean $p$.
    If $p=o(1)$ and $np\to\infty$ as $n\to\infty$, then, with probability $1-o(1)$, $$\Nind(C_{2\ell},G)=(1+o(1))\frac{(np)^{2\ell}}{4\ell}.$$
    For the symmetric algebraic model on either the projective set or its affine chart, fixed $b>a\ge1$ and $m\ge4\ell-1$ therefore give \begin{equation}\label{eq:sharp-induced-cycle-count} \Nind(C_{2\ell},G)=(1+o(1))\frac{q^{2\ell(b-a)}}{4\ell} \end{equation} with probability $1-o(1)$.
\end{theorem}

\begin{proof}
    Put $X=\Nind(C_{2\ell},G)$ and $d=np$.
    Since $|\operatorname{Aut}(C_{2\ell})|=4\ell$, Proposition~\ref{prop:induced-moment-transfer} gives $$\mathbb{E}X=\frac{(n)_{2\ell}}{4\ell} p^{2\ell}(1-p)^{\ell(2\ell-3)} =(1+o(1))\frac{d^{2\ell}}{4\ell},$$ and reduces the second moment to $G(n,p)$.
    Index candidate cycles modulo their automorphisms.
    Vertex-disjoint pairs have independent indicators and contribute $(1+O_\ell(n^{-1}))(\mathbb{E}X)^2$; identical pairs contribute $\mathbb{E}X=o((\mathbb{E}X)^2)$.

    For a distinct overlapping pair, let $r\ge1$ be the number of common vertices and $s$ the number of common cycle edges.
    Conflicting induced requirements give probability zero; otherwise the joint probability is at most $p^{4\ell-s}$.
    Each overlap type has $O_\ell(n^{4\ell-r})$ ordered pairs.
    The common-edge graph consists of paths and isolated vertices, since a cycle component would force the candidates to coincide.
    Thus $s\le r$, and its normalized contribution is $$O_\ell(n^{-r}p^{-s})\le O_\ell((np)^{-r})=o(1).$$
    There are finitely many overlap types, so $\operatorname{Var}(X)=o((\mathbb{E}X)^2)$.
    Chebyshev's inequality proves the result.
    Lemma~\ref{vertex-supported-independence} applies to the algebraic model when $m\ge4\ell-1$, where $n=(1+o(1))q^b$, $p=q^{-a}=o(1)$ and $np\to\infty$.
\end{proof}

Under the same asymptotic assumptions, the argument applies to each fixed $C_s$, $s\ge3$, using independence on at most $2s$ vertices ($m\ge2s-1$ in the algebraic model); we state the even-cycle case to compare it with the following spectral application.

A different conclusion follows from spectral information alone.
The next proposition applies the deterministic supersaturation theorem of D\v{z}avoronok, Gabsdil, Mylet, Popa and Shao \cite[Theorem~1.2]{DGMPSh_induced_cycles}.
It requires an additional density condition and gives a lower bound, rather than an asymptotic count in the sampled model.

\Needspace{13\baselineskip}
\begin{proposition}[A deterministic spectral application]\label{thm:induced-cycles-spectral}
    Fix integers $b>a\ge1$, $\ell\ge2$ and $h\ge1$, and $K>0$.
    Let $G=G_q$ be a sequence of graphs on $n=n(q)=(1+o(1))q^b$ vertices, and put $p=q^{-a}$.
    Suppose that $$\|A(G)-p(J-I)\|\le\beta_h,\qquad \beta_h=K n^{1/2+1/(2h)}p^{1/2}.$$
    If \begin{equation}\label{eq:spectral-induced-condition} (\ell+1)(b-a)>2b+\frac{(\ell-1)b}{h}, \end{equation} then $\Nind(C_{2\ell},G)=\Omega_{\ell,a,b,h,K}((np)^{2\ell})$.
\end{proposition}

\begin{proof}
    Put $d=np$.
    Condition~\eqref{eq:spectral-induced-condition} implies $b/h<b-a$, so $\beta_h=o(d)$ and the average degree is $(1+o(1))d$.
    Choose $\varepsilon=\varepsilon_\ell>0$ from the cited theorem and put $T=\lceil2\beta_h/\varepsilon\rceil$.
    Then $T=\Theta(d^{1/2}n^{1/(2h)})$.
    For $|U|,|W|\ge T$, the spectral estimate gives $$e(U,W)\le p|U||W|+\beta_h\sqrt{|U||W|} \le(p+\beta_h/T)|U||W|\le\varepsilon|U||W|$$ for large $q$.
    Thus the graph is $(1-\varepsilon,T)$-sparse in the terminology of that theorem, with the same convention for overlapping sets.
    The cited theorem requires $$e(G)\ge C_\ell T^{1-1/\ell}n^{1+1/\ell}.$$
    Since $e(G)=(1+o(1))nd/2$, this threshold is met because $$\frac{d}{T^{1-1/\ell}n^{1/\ell}}\longrightarrow\infty \quad\Longleftrightarrow\quad \frac{d^{\ell+1}}{n^{2+(\ell-1)/h}}\longrightarrow\infty,$$ and the latter is exactly \eqref{eq:spectral-induced-condition}.
    The theorem gives $\Nind(C_{2\ell},G)=\Omega_\ell((2e(G)/n)^{2\ell}) =\Omega_\ell(d^{2\ell})$, as required.
\end{proof}

A fixed $h$ satisfying \eqref{eq:spectral-induced-condition} exists precisely when $b/a>(\ell+1)/(\ell-1)$.
This restriction concerns this application of the supersaturation theorem, whereas Theorem~\ref{thm:sharp-induced-cycles} applies for every $b>a$.

\begin{corollary}[Simultaneous realization]\label{cor:induced-cycles-simultaneous}
    Let $(F,R)$ be $(a,b)$-balanced, with $b>a\ge1$, and fix integers $\ell\ge2$ and $h\ge1$.
    Choose an even integer $r>2b/(b-a)$ and a fixed degree $$m\ge\max\{4\ell-1,\,2h-1,\,r-1,\,(b|R|+1)e(F)-1\}.$$
    There are a fixed $t$ and $K>0$ such that, for every sufficiently large $q=2^s$, an $F^{(t)}$-free graph $G_q$ on $\mathbb{F}_q^b$ satisfies $$\begin{aligned} \deg_{G_q}(v)&=(1+o(1))q^{b-a}&&\text{uniformly in }v,\\ \|A(G_q)-q^{-a}(J-I)\|&\le Kq^{(b-a)/2+b/(2h)},\\ \Nind(C_{2\ell},G_q)&=(1+o(1))q^{2\ell(b-a)}/(4\ell). \end{aligned}$$ Here $\ell,h$ are fixed before $m$, and $t$ is chosen only after $m$.
\end{corollary}

\begin{proof}
    Apply Proposition~\ref{prop:symmetric-coefficient-obstruction} with this degree $m$ over $\mathbb{F}_2$ to obtain $t$ and $Q$.
    Then $Q(c)\ne0$ has probability $1-O(q^{-1})$.
    The projective spectral event has probability at least $1/2$ and passes to the affine chart.
    Degree concentration, applied directly on that chart, and the induced-cycle estimate each hold with probability $1-o(1)$.
    The intersection has probability at least $1/2-O(q^{-1})-o(1)>0$.
\end{proof}

Under the full hypotheses of Theorem~\ref{H_r}, fix $\ell$, choose $h>3b/(b-a)$, and then take $m$ satisfying the budget above.
The resulting $t$ and graphs $G_q$ may also be used in the Ramsey argument of Section~5.

\section*{Acknowledgements and disclosure}
The author is grateful to Hong Liu for suggesting the problem and for his guidance and help throughout this research.

ChatGPT assisted with drafting the initial version of this manuscript.
The proof of the extremal input in Theorem~\ref{thm:clean-power-input} was claimed by Epoch AI during benchmarking of GPT-6 Astra and is presented in human-readable form by Conlon et al.~\cite{conlon2026rational}.

\bibliographystyle{abbrv}
\bibliography{Single-graph_rational_exponents_in_mixed_multicolour_Ramsey_theory}
\end{document}